\documentclass[reprint,12pt]{elsarticle}
\usepackage{amssymb}
\usepackage{bm}
\usepackage{ulem}
\usepackage{xcolor}
\usepackage{lineno}

\begin{document}

\begin{frontmatter}

\title{Simple numerical scheme for solving the impregnation equations in a porous pellet}
 
\author{N.V. Peskov}
\corref{*} 
\ead{peskovnick@gmail.com}
\author{T.M. Lysak}
\cortext[*]{Corresponding author}

\address{Faculty of Computational Mathematics and Cybernetics, \\Lomonosov Moscow State University,  Moscow, Russian Federation.}

\begin{abstract}
This paper proposes a numerical scheme for solving a system of convection-reaction-diffusion equations describing the process of preparing a catalyst on a porous support by the impregnation method. In the case of a considered porous spherical pellet, the equations are defined on an interval, one end of which, associated with the front of the impregnating liquid, moves according to a given law. The law of front motion is used to create a consistent space-time grid for discretizing the system. Examples of numerical solutions of the impregnation problem are given, demonstrating the effectiveness of the proposed scheme.
\end{abstract}


\end{frontmatter}

\section{Introduction}
The process of impregnating porous pellets with a solution containing catalytically active elements is a widely used technique for the preparation of supported catalysts \cite{CEC20}. Active elements from the solution are adsorbed onto the pore surface, forming a catalytically active surface of large area. The impregnation process is modelled by an initial boundary value problem for a system of partial differential equations in a domain with a moving boundary \cite{LA85, L01, L08, L22}. 

In the models of supported catalyst preparation using the impregnation method, the issues of liquid imbibition in a porous medium due to capillary action and the impregnation of catalytically active elements are considered separately. The solution of the hydrodynamic problem of liquid imbibition leads to the obtaining of the time-dependent function $r_f(t)$, which defines the position of the advanced liquid front within the pellet. This function depends on the properties of the medium, the model chosen for the capillary pressure of the liquid in the pores and the potential behaviour of the air contained in the pellet.

The impregnation equations, actually the two-phase system of reaction-diffusion equations for the solute concentration in the liquid and the adsorbed species concentration on the pore surface, are defined in the part of the pellet occupied by the advanced liquid. As the inner boundary of the wet part of the pellet moves inwards during impregnation, the system of equations is solved in the changing domain. This circumstance causes definite difficulties for the numerical solution, which can be fixed by introducing a new space coordinate depending on the ratio $r/r_f$. But in fact such a fixing creates other difficulties for the solution of the system.

In this paper we propose a numerical scheme for solving the equations of impregnation in a domain with a moving boundary without fixing the boundary, provided a law of boundary motion is known. The main idea of the scheme is to solve the equations on consistent spatial and temporal grids.

\section{Impregnation equations}
In this paper we are only interested in the algorithm of numerical solution of the impregnation problem. Therefore, the derivation of the equations and the evaluation of the values of the model parameters are not presented. These issues are discussed in detail in \cite{LA85} and \cite{L22}. We take the statement of the problem from \cite{LA85}, preserving the names of the variables and parameters.

In dimensionless form, the impregnation equations for an isotropic spherical pellet can be written as follows:
\begin{eqnarray}
\label{ueq}
&&(1-\rho)^2\frac{\partial}{\partial\tau}(u+\eta\theta) = -Q(\tau)\frac{\partial u}{\partial\rho} + d\frac{\partial}{\partial\rho}\left[ (1-\rho)^2\frac{\partial u}{\partial\rho}\right], \\
\label{theq}
&&\frac{\partial\theta}{\partial\tau} = K^+u\left(1-\theta\right) - K^-\theta.
\end{eqnarray}
Here $\rho$ is the (inverted) dimensionless radial coordinate, $\rho=0$ corresponds to the outer surface of the pellet, $\rho=1$ corresponds to the centre of the pellet; $\tau$ is the dimensionless time. Unknown functions: $u(\rho,\tau0$ is the dimensionless concentration of solute in liquid and $\theta(\rho,\tau)$ is the dimensionless consentration of adsorbed species on pore surface, $0\leq u,\theta \leq 1$. Dimensionless parameters of the model: $\eta$ related to the adsorption capacity of the pore surface, $Q(\tau)$ related to the volumetric rate of liquid imbibition, $d$ related to the diffusivity of the solute in the liquid, $K^+$ and $K^-$ related to the adsorption and desorption rate constants respectively.

The equations (\ref{ueq}), (\ref{theq}) are held in the interval $0<\rho <\rho_f(\tau)$, where $\rho_f(\tau)$ is the position of the liquid front at time $\tau$. $Q$ and $\rho_f$ are related by the equation
\begin{equation}
\label{qrf}
Q(\tau) = (1-\rho_f(\tau))^2\frac{d\rho_f(\tau)}{d\tau}.
\end{equation}
The law of motion of the liquid front $\rho_f(\tau)$ can be obtained from the solution of the problem of fluid flow in a porous medium under the action of capillary forces.

The solution to the system of equations (\ref{ueq}), (\ref{theq}) is contingent upon the provision of initial and boundary conditions.
Boundary conditions: at $\rho=0$ -- the Dancwerst conditions
\begin{equation}
\label{bc0}
d\frac{\partial u}{\partial \rho} = Q(\tau)(u-u^0),
\end{equation}
where $u^0=1$ is the solute concentration in the  liquid in which the pellet is immersed. 
At $\rho = \rho_f(\tau)$ the Neumann conditios
\begin{equation}
\label{bc1}
\frac{\partial u}{\partial \rho} = 0.
\end{equation}

At $\tau=0$ the pellet is only immersed in the solution, impregnation is only beginning. It is assumed that $\rho_f(0)=0$ and the system (\ref{ueq}), (\ref{theq}) is not yet defined. Therefore, there is no place to set the initial conditions. In the numerical solution, we bypass this difficulty as it is done in \cite{LA85}.

\section{Numerical scheme}

Our proposed numerical solution of the system (\ref{ueq}), (\ref{theq}) is defined on the coordinated grids along the radius of the pellet and along time. These grids are constructed as follows.

Let as a result of the solution of the hydrodynamic problem obtained (analytically or numerically) function $\rho_f(\tau)$, the time dependence of the position of the liquid front in the pellet. Suppose that $\rho_f(\tau)$ is monotonically increasing function. The path of the front on the plane ($\tau, \rho$), -- curve ($\tau, \rho_f(\tau)$) starts at the point (0,0) and ends at some point $(\tau_e, \rho_e)$. Let $L$ be the length of the path, we divide it into $N$ parts of equal length $L/N$ by points $(\tau_i, \rho_i=\rho_f(\tau_i))$. $i = 0, ..., N$. (If the terminal point $\rho_e$ is reached asymptotically at $\tau \to \infty$, it can be excluded from consideration.) The time moments $\tau_i$, $i=1,2,...,N$ make up the time grid. At each $i$, the points $\bar{\rho}_j = (\rho_{j-1}+\rho_j)/2$, $j = 1,2,..., i$ make up a radial grid in the wet region.

Using the finite difference method, we write the system (\ref{ueq}), (\ref{theq}) in discrete form on tht non-uniform grid ($\tau_i, \bar{\rho}_j$). The implicit Euler scheme is used to approximate the time derivatives. To calculate the radius derivatives at the extreme points of the grid, we use the boundary conditions (\ref{bc0}), (\ref{bc1}). Let the solution of the difference equations for the moment of time $\tau_{i-1}$ be known. Since equation (\ref{theq}) is nonlinear, we will look for the solution of the system for the next time instant $\tau_i$ by the iteration method.

In each $k$-th iteration, the descrete form of the equation (\ref{theq}) for $\theta^{(k)}$ is solved first, and the value $u^{(k-1)}$ from the previous iteration is substituted in place of $u$, or, for the first iteration, the value $u^{(i-1)}$ from the previous time step. The obtained values of $\theta^{(k)}$ are substituted into the descrete form of the equation (\ref{ueq}) and $u^{(k)}$ is calculated. Iterations are stopped when the norm of the difference between successive iterations becomes less than a given value.

For the described scheme to work, it is necessary to set the values on the first time plane to $\tau = \tau_1$. On this layer there is only one point $\bar{\rho}_1$ in the radial grid, so we need to set $u_1^1$ and $\theta_1^1$. In impregnation models, the fluid velocity is very high at the beginning of the process and therefore $\tau_1$ is very small. As in \cite{LA85}, we assume that for the time $\tau_1$ no significant amount of substance has had time to adsorb onto the pore surface and assume $\theta_1^1 = 0$, then $u_1^1 = 1$. These values are used as initial data for solving the scheme.

\section{Numerical example}

To demonstrate the operation of the scheme, we have chosen the apparently simplest model of hydrodynamics of impregnation, in which the capillary pressure $P_c$ of the liquid in the pores is considered constant and equal to $P_c = \sigma P_a$, where $P_a$ is the atmospheric pressure, $\sigma$ is the coefficient, and the air in the pores of the pellet does not resist the movement of the liquid.
\begin{figure}[h]
\centering
\includegraphics{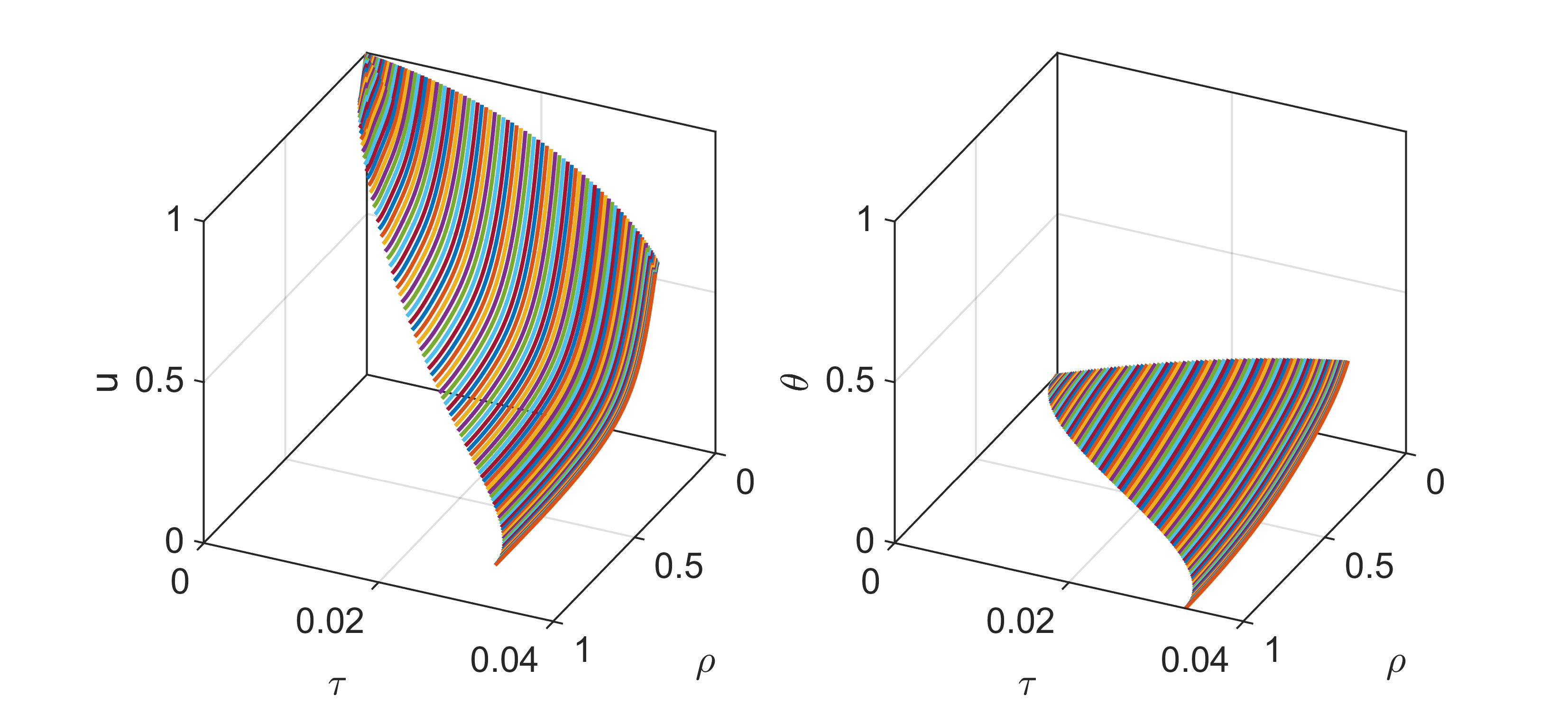}
\caption{Solution of system (1), (2) for $K^+=10$.}
\label{fg10}
\end{figure}
\begin{figure}[h]
\centering
\includegraphics{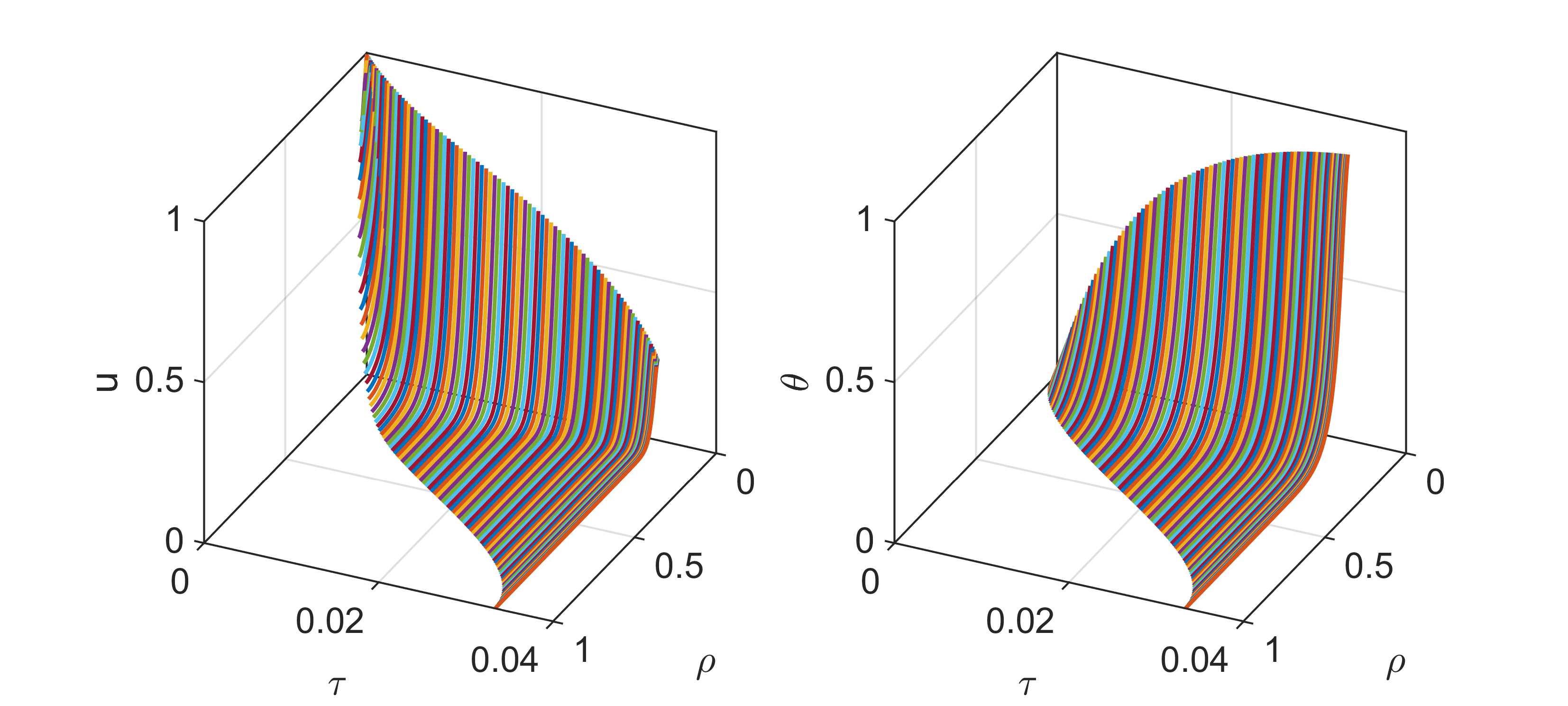}
\caption{Solution of system (1), (2) for $K^+=100$.}
\label{fg100}
\end{figure}
\clearpage

In such a case, the liquid reaches the centre of the pellet in a  time $\tau_e = 1/(6\sigma)$ and the position of the liquid front $\rho_f$ is related to time $\tau$ by the equation 
\begin{equation}
\tau = \frac{1}{6\sigma}\left(3\rho_f^2 - 2\rho_f^3\right).
\label{plrf}
\end{equation}
and from Eq. (\ref{qrf}) one obtains
\begin{equation}
Q(\tau) = \sigma\frac{1-\rho_f}{\rho_f}.
\label{Q}
\end{equation}

The following examples of solutions are obtained with the following parameter values. $\sigma = 5$, so $\tau_e=1/30$. $\eta =6$, $d=0.1$. $K^-=0.1$, and $N=1000$.
 The examples differ in the adsorption constant.  In the first example, Fig. 1, $K^+ = 10$, in the second example. Fig. 2, the adsorption constant is 10 times larger, $K^+ = 100$. 
 
A simple test for correctness of a numerical solution has been proposed in \cite{LA85}. It compares two quantities, the amount of solute $M_1(\tau)$ entered in the pellet for time $\tau$:
\[M_1(\tau) = u^0\int_0^\tau{Q(\xi)\,d\xi},\]
and the amount of of solute $M_2(\tau)$ dissolved in liquid and adsoebed on pore surface in the pellet at time $\tau$,
\[M_2(\tau)=\int_0^{\rho_f(\tau)}{[u(\xi,\tau)+\eta\theta(\xi,\tau)]\,(1-\xi)^2d\xi}\].
For the correct solution it must be $M_1(\tau)=M_2(\tau)$ for each $\tau$.

The next Fig, 3 shows the graphs of $M_1$ ans $M_2$ for the second numerical example.
\begin{figure}[h]
\centering
\includegraphics{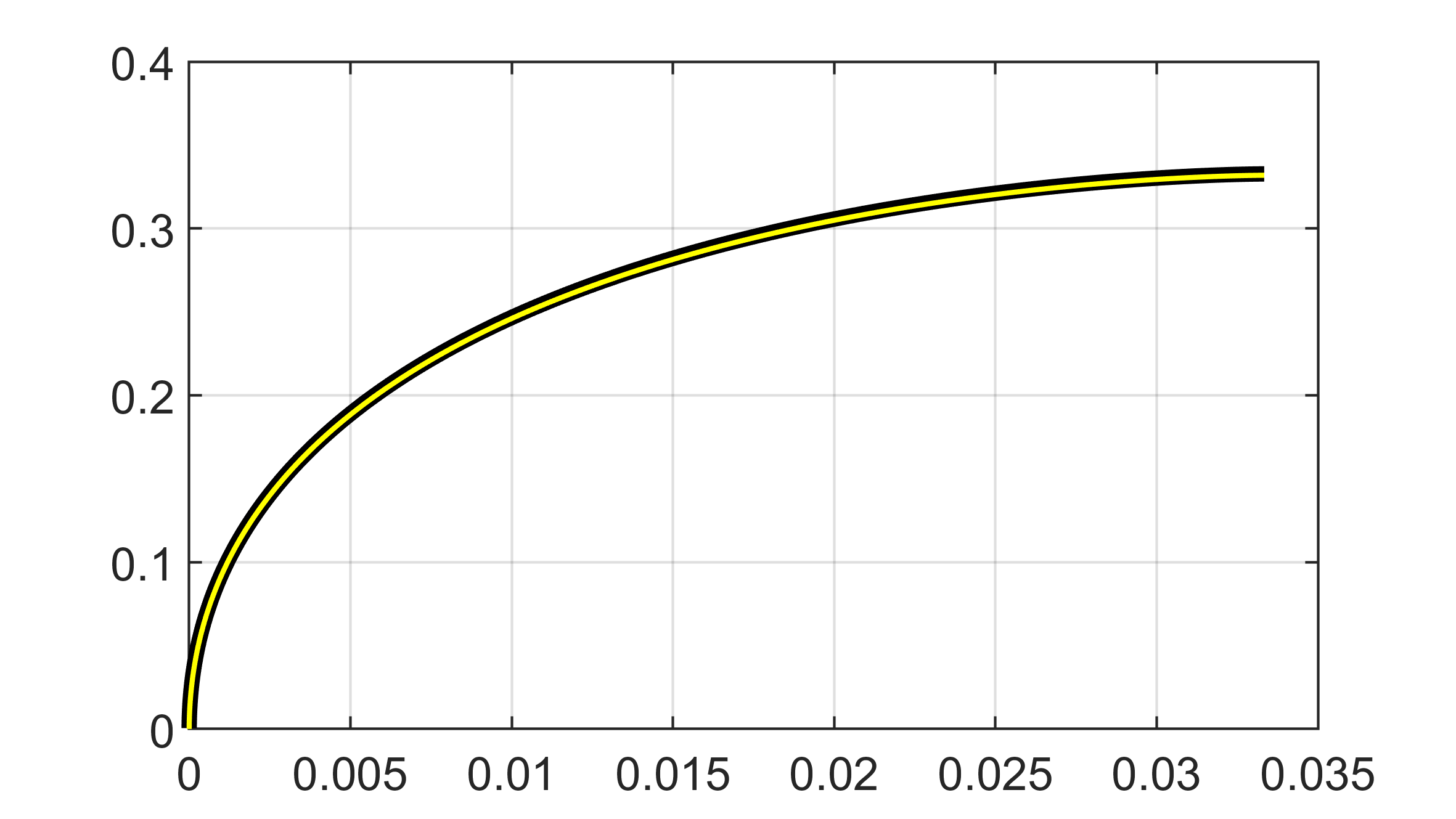}
\caption{$M_1$ -- thick black line, $M_2$ -- thin yellow line.}
\end{figure}

Numerical examples were solved in MATLAB. Iterations were stopped if the uniform norm of the difference of the solution from two consecutive iterations became less than 10$^{-6}$. Solving the examples on an ordinary PC takes a few seconds.

\section{Concluding remarcs}

We presented a numerical scheme for solving the initial boundary value problem for an convection-reaction-diffusion type equations (1), (2) defined on the interval $(0, a(\tau))$, where $a(\tau)$ is a given function of time $\tau$. It is assumed that $a(\tau)$ is a monotonically increasing function without singularities. The finite-difference version of the differential equations is written on a consistent space-time grid $(\tau_i, \bar{\rho}_j)$. $i = 1, 2, ..., N$, $j = 1,2,..., i$. The nodes of the grid are defined by dividing the curve $9\tau, a(\tau)$ into $N$ segments. In the given examples, by dividing the graph into $N$ segments of equal length. The number $N$ and the way of curve division are determined by the required accuracy of the solution and stability of the scheme on the interval $(0,\tau_e)$.

\end{document}